\documentclass[12pt]{amsart}
\usepackage{amssymb}
\usepackage{graphicx}

\numberwithin{equation}{section}

\renewcommand{\P}{{\mathbb P}}
\newcommand{\R}{{\mathbb R}}
\newcommand{\RP}{\R\P}

\newcommand{\QED}{\hfill\raisebox{-5pt}{\includegraphics[height=14pt]{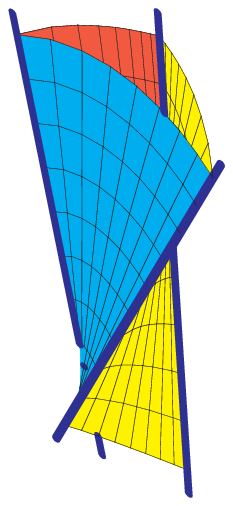}}\vspace{6pt}}

\newtheorem{thm}{Theorem}
\newtheorem{cor}[thm]{Corollary}
\newtheorem{prop}[thm]{Proposition}

\newtheorem{remark}[thm]{Remark}
\newenvironment{rem}{\begin{remark}\rm}{\end{remark}}

\newtheorem{example}[thm]{Example}
\newenvironment{ex}{\begin{example}\rm}{\end{example}}

\newtheorem{definition}[thm]{Definition}
\newenvironment{defn}{\begin{definition}\rm}{\end{definition}}

\title[Cremona Convexity and Santal\'o's Theorem]{Cremona Convexity, Frame
Convexity,\\ and a Theorem of Santal\'o}

\author[J.~E.~Goodman]{Jacob E. Goodman}
\address{Department of Mathematics\\
         City College (CUNY)\\
         New York, NY 10031\\
         USA}
\email{jegcc@cunyvm.cuny.edu}
\urladdr{http://www.ccny.cuny.edu/mathematics/faculty/faculty.htm}
\author[A.~Holmsen]{Andreas Holmsen}
\address{Matematisk institutt\\
         Universitetet i Bergen\\
         Johannes Bruns gt.\ 12\\
         5008 Bergen\\
         Norway}
\email{andreash@mi.uib.no}
\urladdr{http://www.mi.uib.no/\~{}andreash/index.html}
\author[R.~Pollack]{Richard Pollack}
\address{Department of Mathematics\\
         Courant Institute of Mathematical Sciences\\
         New York University\\
         251 Mercer St.\\
         New York, NY 10012\\
         USA}
\email{pollack@cims.nyu.edu}
\urladdr{http://www.math.nyu.edu/faculty/pollack/}
\author[K.~Ranestad]{Kristian Ranestad}
\address{Matematisk institutt\\
         Universitetet i Oslo\\
         PO Box 1053, Blindern\\
         NO-0316 Oslo\\
         Norway}
\email{ranestad@math.uio.no}
\urladdr{http://www.math.uio.no/\~{}ranestad}
\author[F.~Sottile]{Frank Sottile}
\address{Department of Mathematics\\
         Texas A\&M University\\
         College Station\\
         TX \ 77843\\
         USA}
\email{sottile@math.tamu.edu}
\urladdr{http://www.math.tamu.edu/\~{}sottile}

\thanks{Research of first author supported in part by NSA grant 
MDA904-03-I-0087 and PSC-CUNY grant 65440-0034;
research of second author 
supported in part by the Mathematical Sciences
Research Institute, Berkeley;
research of third author supported in part by NSF grant CCR-9732101;
research of fifth author supported in part by NSF CAREER
  grant DMS-0134860 and the Clay Mathematical Institute}
\begin{document}

\begin{abstract}
 In 1940, Luis Santal\'o proved a Helly-type theorem for line transversals
 to boxes in $\R^d$.
 An analysis of his proof reveals a convexity structure for ascending lines
 in $\R^d$ that is isomorphic to the ordinary notion of convexity in a convex
 subset of $\R^{2d-2}$.
 This isomorphism is through a Cremona transformation on the Grassmannian of
 lines in $\P^d$, which enables a precise description of the convex hull
 and affine span of up to $d$ ascending lines:
 the lines in such an affine span turn out to be the rulings of certain classical
 determinantal varieties.
 Finally, we relate Cremona convexity to a new convexity structure that we
 call {\em frame convexity}, which extends to arbitrary-dimensional flats
 in $\R^d$.
\end{abstract}

\maketitle

\section*{Introduction}
A cornerstone of the theory of convexity is Helly's theorem~\cite{He1923}, which
states that if there is a point common to every $d{+}1$ (or fewer) members of a
collection of compact convex sets in $\R^d$, then the whole family has a point
in common.
Vincensini~\cite{Vi1935} asked if there was a similar result
for line secants to compact convex sets in the plane.
In 1940, Santal\'o showed that this was impossible in general, found
a sufficient condition for the result to hold for lines in the plane, and 
generalized that to
lines in $\R^d$, giving a Helly-type theorem for line transversals to
boxes~\cite{Sa1940}.
A {\it box} in $\R^d$ is a product of $d$ intervals, one in each coordinate
direction.

\begin{thm}[Santal\'o]\label{T:Santalo_Lines}
 If every $2^{d-1}(2d{-}1)$ (or fewer) members of a collection of boxes in
 $\R^d$ admit a line transversal, then there is a line transversal to all the
 boxes.
\end{thm}

This was the first of a great many results in what is now known as
geometric transversal theory~\cite{DGK63,Ec93,GPW93,We99,We04}.
Santal\'o's elementary proof was based on an idea that Radon had used to establish
Helly's theorem~\cite{Ra1921}.
Gr\"unbaum~\cite{Gr58} showed that these numbers are the best possible;
there are 6 squares in the plane, every five of which have a line
transversal, but all six do not.
If we ignore the factor $2^{d-1}$, 
which is simply the number of direction classes for lines in $\R^d$ that are
not parallel to any coordinate hyperplane,
the number in Santal\'o's theorem is the
Helly number for dimension $2d{-}2$, which is the dimension of the space of
lines in $\R^d$.
This is not a coincidence---for as we shall see, Santal\'o's theorem is
essentially Helly's theorem in disguise.

A line in $\R^d$ is {\it weakly ascending} if all non-zero coordinates of its
direction vector have the same sign.
If no coordinates of its direction vector vanish, then such a line
is {\em ascending}.
A {\em frame-convex} set of (ascending) lines is the set of lines that
meet every box in a given collection of boxes.  (Since the notion of a `box'
is dependent on the choice of coordinate frame, this notion of convexity
is as well.)
This is related to a general theory of convexity on Grassmannians
due to Goodman and Pollack~\cite{GP95}
in which a set $\mathcal{K}$ of $k$-flats is convex if whenever a $k$-flat 
$K$ meets every convex body meeting every member of $\mathcal{K}$, then
$K$ lies in $\mathcal{K}$.
For frame convexity, we restrict the convex bodies to be boxes and
work with lines rather than $k$-flats.
Weakly ascending lines admit natural Cremona coordinates (described in
Section~\ref{S:Santalo}) which identify them with a convex set in $\R^{2d-2}$.

\begin{thm}\label{T:Convexity}
 The Cremona coordinates of a frame-convex set of weakly ascending lines form a 
 convex set.
\end{thm}

This elementary result can be seen in a close reading of Santal\'o's arguments.
We give a short proof of Theorem~\ref{T:Convexity} and also of Santal\'o's
theorem in Section~\ref{S:Santalo}.
The converse of Theorem~\ref{T:Convexity} is not true in general.
A set $S$ of weakly ascending lines is {\it Cremona-convex} if the set of Cremona
coordinates of its members form a convex set.
Suppose that $\ell$ and $\ell'$ are ascending lines in $\R^d$ that
meet in a point.
Then their Cremona-convex hull is 1-dimensional, but their frame-convex hull may
have any dimension between 1 and $d{-}1$.
The case $d=3$ is described in Example~\ref{Ex:lines_meet}.
Despite this, these two notions of convexity coincide for ascending skew lines
that satisfies the additional condition that the images of the direction 
vectors by the orthogonal projection to any coordinate plane
are distinct 
(see Corollary~\ref{cor:general}).

For ascending lines, the Cremona coordinates of Theorem~\ref{T:Convexity} come
from a particular Cremona transformation on Pl\"ucker space that transforms the
Grassmannian of lines  into a linear space.
This linearizing Cremona transformation is
related to Kapranov's identification of the Chow quotient of the
Grassmannian of lines as the space $M_{0,n}$ of $n$ marked points on
$\P^1$~\cite{Ka93}, and it underlies some recent work on the tropical
Grassmannian of lines~\cite{AK03,SpSt}.

\begin{figure}[htb]
 \[
  \includegraphics[width=7cm]{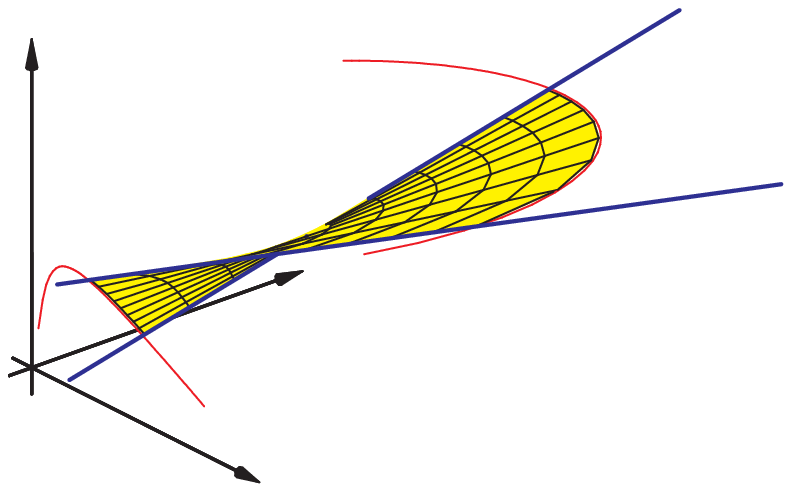}\qquad
   \includegraphics[width=7cm]{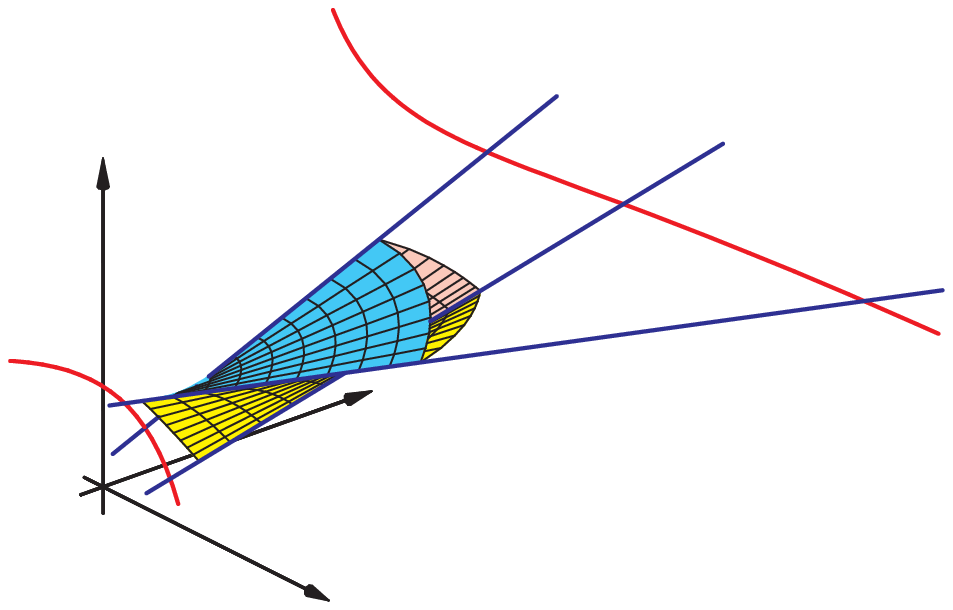}
 \]
  \caption{Cremona-convex hulls}\label{F:Cremona_Interval}
\end{figure}

We obtain
a precise geometric description of the
Cremona-convex hull and Cremona-affine span of a set of lines.
For example, consider two skew ascending lines in $\R^3$ such that there is a unique line parallel to 
each axis meeting the two given lines.
Then these three axis-parallel lines lie on a unique doubly ruled
quadric, either a hyperboloid of one sheet or a hyperbolic paraboloid.
The three axis-parallel lines lie in one ruling and the two original
lines lie in the other.
The two original lines determine two intervals in their ruling
(three in the case of a hyperbolic paraboloid). 
One of those intervals consists solely of ascending lines, and this interval is
the frame-convex interval between the original two lines.
This is shown in the figure on the left in Figure~\ref{F:Cremona_Interval}.
The Cremona-affine span of three general lines in $\R^3$ consists of the secant
lines to a particular rational cubic curve.
This is illustrated in the figure on the right in
Figure~\ref{F:Cremona_Interval}, where we display two branches of the
rational cubic and lines that are the extreme points of the
Cremona-convex hull of three lines. 
In Section~\ref{S:Cremona-affine} we explain this more generally in $\R^d$. 

Replacing boxes by parallelepipeds with parallel edges gives the coordinate-free
version of frame convexity. 
This was the original context of Santal\'o's results.
Santal\'o also proved a Helly-type theorem for hyperplane transversals, which
follows from the observation that for ascending hyperplanes
(suitably defined), frame convexity coincides with a natural notion of
convexity given by the coefficients of the defining equation.

Santal\'o asked if his results could be extended to $k$-flats in $\R^d$, for
$1<k<d{-}1$.
We do not know if this is possible, but feel that it is unlikely.
In part that is because the linearizing Cremona transformation does not generalize 
to the Grassmannian with $1<k<d-2$.
Also, we show that the frame-convex hull of two general 2-flats in $\R^4$ or
$\R^5$ (suitably defined) is just the two original 2-flats, and thus frame
convexity in this context has a completely different character than for lines or
hyperplanes.  
On the other hand, since the set of $(d{-}2)$-flats in $\P^d$ is isomorphic to
the set of lines in $\P^d$, the notion of frame convexity for ascending lines
can be transferred to $(d{-}2)$-flats and gives it a convexity structure.
We describe this in Section~\ref{S:Frame}, and derive the corresponding
theorem of Helly type.

We thank Bernd Sturmfels, who pointed out to us the relation of the Cremona
transformation to Chow quotients and to the tropical Grassmannian.

\section{Convexity for weakly ascending lines in $\R^d$}\label{S:Santalo}

A weakly ascending line $\ell$ in $\R^d$ is determined by a point $x\in\ell$ and
its direction vector $v$, which has non-negative coordinates.
Call such a pair $(x,v)$ {\it rectilinear coordinates} for $\ell$.
The {\it Cremona coordinates} $(y,w)$ for a line $\ell$ are defined as follows
 \begin{equation}\label{E:Cremona_Coords}
   \left\{
      \begin{array}{rclcrclcl}
       y_i&=&{\displaystyle \frac{x_i}{v_i}}&{\textrm \ and\ }&
       w_i&=&{\displaystyle \frac{1}{v_i}}&
     \quad\mbox{if}&v_i\neq 0\\
       y_i&=&x_i&{\textrm \ and\ }&w_i&=&v_i\ =\ 0&
     \quad\mbox{if}&v_i= 0 \rule{0pt}{15pt}
     \end{array}\right.  \qquad{\textrm for}\quad i=1,\dotsc,d\,.
 \end{equation}
The first Cremona coordinate $y$ of a weakly ascending line is well-defined
modulo translation by the vector $(1,\dotsc,1)$ and the second Cremona
coordinate $w$ is well-defined modulo multiplication by a positive scalar.
If we remove these ambiguities by requiring $\sum_i y_i=0$ and $\sum_i w_i=1$, 
then the Cremona coordinates take values in a convex subset of a 
$(2d{-}2)$-dimensional affine subspace of $\R^{2d}$.

\begin{defn}
  A set of weakly ascending lines is {\it Cremona-convex} if its set of Cremona
  coordinates is convex.
\end{defn}

 These definitions imply that the dimension of the Cremona-convex
 hull of $k$ lines in $\R^d$ is at most $k{-}1$, with equality only if the Cremona
 coordinates are affinely independent.
 With this definition, Theorem~\ref{T:Convexity} becomes the following.\medskip

\noindent{\bf Theorem~\ref{T:Convexity}$'$.}
 {\it
  A frame-convex set of weakly ascending lines is Cremona-convex.
}\medskip

 Since an ascending line transversal to a collection of boxes is a common
 member of each of the frame-convex sets of ascending lines meeting each box,
 Helly's theorem applied to such convex sets gives a local version of
 Santal\'o's theorem.

\begin{cor}\label{C:Local_Santalo}
 If every $2d{-}1$ (or fewer) members of a collection of boxes in $\R^d$ admit
 an ascending line transversal, then there is an ascending line transversal to
 all the boxes.
\end{cor}

\noindent{\it Proof of Theorem~$2$.}
 A box in $\R^d$ is the coordinatewise interval between its (coordinatewise)
 minimal and maximal points.
 The rectilinear coordinates $(x,v)$ of ascending lines that meet a box
 with minimal point $a$ and maximal point $b$ satisfy the inequalities
 \begin{equation}\label{E:ineq}
  \begin{array}{rcl}
  {\displaystyle \max_{i\colon v_i\neq 0}\left\{\frac{a_i-x_i}{v_i}\right\}}
   &\leq&
  {\displaystyle \min_{i\colon v_i\neq 0}\left\{\frac{b_i-x_i}{v_i}\right\}}\,,
     \qquad\mbox{and}\\
    a_i\ \leq& x_i&\leq\ b_i \qquad \mbox{if}\quad v_i=0\,.\rule{0pt}{15pt}
  \end{array}
 \end{equation}
 Indeed, a line with rectilinear coordinates $(x,v)$ meets the box
 if and only if the system of inequalities
\[
  a_i\ \leq\ x_i+t\cdot v_i\ \leq\ b_i\,\qquad
   \textrm{for}\quad i=1,\dotsc,n
\]
 has a solution in $t$, for then $x+t\cdot v$ lies in the box.
 The feasibility of this system for weakly ascending lines is expressed
 by~\eqref{E:ineq}. 

 In terms of Cremona coordinates, the inequalities~\eqref{E:ineq} become
\[
  \begin{array}{rcl}
  {\displaystyle \max_{i\colon w_i\neq 0}\{a_iw_i-y_i\}}
   &\leq&
 {\displaystyle  \min_{i\colon w_i\neq 0}\{b_iw_i-y_i\}}\,,
     \qquad\mbox{and}\\
    a_i\ \leq& y_i&\leq\ b_i \qquad \mbox{if}\quad w_i=0\,,\rule{0pt}{15pt}
  \end{array}
 \]
 which define a convex set.
\QED

\begin{rem}\label{R:Santalo_proof}
 Cremona coordinates of a line are the coefficients of some special equations
 for this line.
 (They appear in Santal\'o's proof in this form.)
 Let $z_1,\dotsc,z_d$ be coordinates for $\R^d$, and consider a line $\ell$ with
 rectilinear coordinates $(x,v)$ where no coordinate of $v$ vanishes.
 Scaling if necessary, we may assume that $v_1=1$,
 and then transforming $x$ by $v$, that $x_1=0$.
 Then the line $\ell$ given by $x+tv$ is cut out by the linear forms
\[
   v_i\cdot z_1\ -\ z_i\ +\ x_i\ =\ 0\,,\qquad\mbox{for}\quad i=2,\dotsc,d\,.
\]
 If we divide by $v_i$, this becomes
 \begin{equation}\label{E:Cremona_revelaled}
    z_1\ -\ w_i\cdot z_i\ +\ y_i\ =\ 0\quad\mbox{if } w_i\neq 0
     \qquad\mbox{for }i=2,\dotsc,d\,.
 \end{equation}
 Santal\'o's proof uses a convex combination of the coefficients $(b_i,c_i)$
 of the linear equations of a line that have the form
\[
  a\cdot z_1\ +\ b_i\cdot z_i\ =\ c_i\qquad i=2,\dotsc,d\,,
\]
 which we recognize from~\eqref{E:Cremona_revelaled} as essentially the Cremona
 coordinates.

 While Cremona coordinates of a line have this simple natural explanation as the
 coefficients of linear forms defining the line, that ad hoc realization
 belies their naturality as the coordinates for the Grassmannian in the
 linearizing Cremona transformation given in Section~\ref{S:linearizing}.
 Using this Cremona transformation, we will obtain a precise description of
 the Cremona-convex hull of a finite set of ascending lines.
\end{rem}

\noindent{\it Proof of Santal\'o's theorem}.
 The closed orthants in the hyperplane at infinity are indexed by the
 signs of coordinates of points in their interiors.
 These vectors $\varepsilon\in\{\pm1\}^d$ are defined modulo
 multiplication by $\pm 1$.
 A line $\ell$ has {\it sign} $\varepsilon$ if it meets the hyperplane at
infinity
 in the closed orthant corresponding to $\varepsilon$.
 Equivalently, if the numbers
\[
   \varepsilon_i v_i\qquad i=1,\dotsc,d
\]
 are either all non-negative or all non-positive, where $v$ is the direction
 vector of $\ell$.

 Our previous notions and results hold, {\it mutatis mutandis}, for lines
 with a fixed sign.
 Every line has Cremona coordinates~\eqref{E:Cremona_Coords}, and we call
 a set of lines with the same sign Cremona-convex if its corresponding
 set of Cremona coordinates is convex.
 The full set of lines with a given sign is Cremona-convex.\medskip

\noindent{\bf Corollary~\ref{C:Local_Santalo}$'$.}
{\it
  If every $2d{-}1$ (or fewer) members of a collection of boxes in $\R^d$ admit
  a line transversal with sign $\varepsilon$, then the collection has a line
  transversal with sign $\varepsilon$.\medskip
}

 Suppose that we have a family of boxes
 such that every subset $S$ of size at most $2^{d-1}(2d{-}1)$ admits a line
 transversal.
 Then there is a sign $\varepsilon$ such that every collection $S$ of
 $2d{-}1$ boxes from our family admits a line transversal with sign
 $\varepsilon$.
 Suppose this were not the case.
 Then for every sign $\varepsilon$, there is a set $S_\varepsilon$ of
 $2d{-}1$ boxes that does not admit a line transversal with sign
 $\varepsilon$.
 But then the union of the sets $S_\varepsilon$ cannot admit a line transversal,
 which is a contradiction, as this union contains at most $2^{d-1}(2d{-}1)$
 boxes.
 Santal\'o's theorem now follows by its local version,
 Corollary~\ref{C:Local_Santalo}$'$.
\QED

\begin{ex}\label{Ex:lines_meet}
 Suppose that $\ell\neq\ell'$ are ascending lines in $\R^3$ that meet
 in a point, which we assume is the origin.
 Let $(a,b,c)$ be the direction vector of $\ell$ and $(a',b',c')$ that for $\ell'$.
 Then Cremona coordinates for $\ell$ and $\ell'$ are 
\[
   \left[(0,0,0),\ \left(\frac{1}{a},\frac{1}{b},\frac{1}{c}\right)\right]
   \qquad\mbox{\rm and}\qquad
   \left[(0,0,0),\ 
      \left(\frac{1}{a'},\frac{1}{b'},\frac{1}{c'}\right)\right]\,,
\]
 and their Cremona-convex hull consists of lines through the origin with
 direction vector
\[
  \left(
    \frac{1}{\frac{p}{a}+\frac{p'}{a'}},\ 
    \frac{1}{\frac{p}{b}+\frac{p'}{b'}},\ 
    \frac{1}{\frac{p}{c}+\frac{p'}{c'}}\right)\ ,
   \qquad\mbox{\rm where}\quad p,p'\geq 0\quad\mbox{\rm and}\quad p+p'=1\,.
\]
These direction vectors $(x,y,z)$ satisfy the quadratic equation
 \begin{equation}\label{Eq:Cremona-hull}
    \left|\begin{matrix} 
     xy&xz&yz\\bc&ac&ab\\b'c'&a'c'&a'b'\end{matrix}\right|\ =\ 0\,,
 \end{equation}
where $|\cdot|$ denotes determinant.

This defines is an irreducible quadric unless one of the $2\times 2$-minors involving
the last two rows of the matrix~\eqref{Eq:Cremona-hull} vanishes.
When that happens, $\ell$ and $\ell'$ span a plane that contains a coordinate
axis. 
Otherwise, the Cremona-convex hull of $\ell$ and $\ell'$ consists of some lines
that rule the cone with apex the origin that is defined by the
equation~\eqref{Eq:Cremona-hull}. 
Specifically, the lines $\ell$ and $\ell'$ define two intervals in the ruling of
this cone, one of which consists entirely of ascending lines, and this
interval is their Cremona-convex hull.

The situation is different for the frame-convex hull of $\ell$ and $\ell'$.

\begin{prop}\label{P:frame2big}
 Suppose that 
 \begin{equation}\label{Eq:sign}
   \left|\begin{matrix}a&b\\a'&b'\end{matrix}\right|\cdot
   \left|\begin{matrix}a&c\\a'&c'\end{matrix}\right|\ \leq\ 0\,.
 \end{equation}
 Then the frame-convex hull of $\ell$ and $\ell'$ is the set of ascending
 lines through the origin with direction vector $(\alpha,\beta,\gamma)$
 satisfying 
 \begin{equation}\label{Eq:frame2big}
   \frac{b'}{a'}\leq \frac{\beta}{\alpha}\leq\frac{b}{a}
   \qquad\mbox{\rm and}\qquad
   \frac{c}{a}\leq \frac{\gamma}{\alpha}\leq\frac{c'}{a}\,.
 \end{equation}
 This is a $2$-dimensional family of lines, except when one of the
 determinants~\eqref{Eq:sign} vanishes, and then it is one-dimensional.
 Since $\ell\neq\ell'$, at most one determinant will vanish.
\end{prop}

\begin{rem}
 The assumption of Proposition~\ref{P:frame2big} can always be satisfied 
 by permuting the coordinate directions.
 Indeed, we have that 
\[
  \mbox{\rm either}\qquad 
   \left|\begin{matrix}a&b\\a'&b'\end{matrix}\right|\cdot
   \left|\begin{matrix}a&c\\a'&c'\end{matrix}\right|\ \leq\ 0
  \qquad \mbox{\rm or}\qquad 
   \left|\begin{matrix}a&b\\a'&b'\end{matrix}\right|\cdot
   \left|\begin{matrix}a&c\\a'&c'\end{matrix}\right|\ >\ 0\ .
\]
 In the second case, 
\[
  \left( \left|\begin{matrix}b&a\\b'&a'\end{matrix}\right|\cdot
         \left|\begin{matrix}b&c\\b'&c'\end{matrix}\right|\right) \cdot
  \left( \left|\begin{matrix}c&b\\c'&b'\end{matrix}\right|\cdot
         \left|\begin{matrix}c&a\\c'&a'\end{matrix}\right|\right)
  \ =\  
 \left( \left|\begin{matrix}a&b\\a'&b'\end{matrix}\right|\cdot
        \left|\begin{matrix}a&c\\a'&c'\end{matrix}\right|\right) \cdot
 \left( \left|\begin{matrix}b&c\\b'&c'\end{matrix}\right|\cdot
        \left|\begin{matrix}c&b\\c'&b'\end{matrix}\right|\right)
   \ \leq\ 0\ .
\]
 Thus, one of the three products is non-positive
\[
   \left|\begin{matrix}a&b\\a'&b'\end{matrix}\right|\cdot
   \left|\begin{matrix}a&c\\a'&c'\end{matrix}\right|\ ,\qquad
   \left|\begin{matrix}b&a\\b'&a'\end{matrix}\right|\cdot
   \left|\begin{matrix}b&c\\b'&c'\end{matrix}\right|\ ,\qquad
   \left|\begin{matrix}c&b\\c'&b'\end{matrix}\right|\cdot
   \left|\begin{matrix}c&a\\c'&a'\end{matrix}\right|\ .
\]
\end{rem}

\noindent{\em Proof of Proposition~$\ref{P:frame2big}$.}
 Interchanging the last two coordinates if necessary, we may assume from 
 ~\eqref{Eq:sign} that $ab'-a'b\leq 0\leq ac'-a'c$, and so
\[
  \frac{b'}{a'}\leq \frac{b}{a}\qquad\mbox{\rm and}\qquad
  \frac{c}{a}\leq \frac{c'}{a'}\,.
\]
 Let $B_0$ be the degenerate box with minimum and maximum points
 \begin{equation}\label{Eq:other_points}
  \left(1,\frac{b'}{a'},\frac{c}{a}\right)\qquad\mbox{\rm and}\qquad
  \left(1,\frac{b}{a},\frac{c'}{a'}\right)\ .
 \end{equation}
 The set of lines through the origin that meet $B_0$ are exactly those whose 
 whose direction vectors satisfy~\eqref{Eq:frame2big}. 
 In particular, $B_0$ meets $\ell$ and $\ell'$.
 The origin is a (degenerate) box meeting $\ell$ and $\ell'$, so their
 frame-convex hull consists of lines through the origin that meet
 $B_0$, as well as every other box meeting both $\ell$ and $\ell'$.

 The ascending lines through the origin that meet a given box $B$ form a
 (double) convex cone with apex the origin.
 We show that if $\ell$ and $\ell'$ meet $B$, then this cone contains the cone
 over $B_0$, which will complete the proof.
 Scaling coordinates and reflecting $B$ in the origin if necessary, we may
 assume that $B$ contains the points
\[
  \left(s, \frac{sb}{a}, \frac{sc}{a}\right)\ \in\ \ell
   \qquad\mbox{\rm and}\qquad
  \left(1, \frac{b'}{a'}, \frac{c'}{a'}\right)\ \in\ \ell'\,,
\]
 for some $s>0$.
 
 Interchanging $\ell$ and $\ell'$ if necessary, we may assume that $s\geq 1$.
 Then the point $(s, \frac{sb'}{a'}, \frac{sc}{a})$ lies on the edge of $B$ 
 beween its vertices $(s, \frac{b'}{a'}, \frac{sc}{a})$ and 
 $(s, \frac{sb}{a}, \frac{sc}{a})$.
 Thus the line $m$ through the first point of~\eqref{Eq:other_points} meets $B$.
 In a similar fashion, the line $m'$ through the second point
 of~\eqref{Eq:other_points} also meets $B$. 
 These four lines $\ell,\ell',m$, and $m'$ generate the cone over $B_0$ with
 apex the origin, which completes the proof.
 \QED
\end{ex}

 Thus the frame-convex hull of two general ascending lines that meet in a point $p$
 fills out a 2-dimensional quadrilateral cone with vertex $p$. 
 This degenerates to a 1-dimensional cone if the affine span of the 
 lines contains a line through $p$ that is parallel to a coordinate axis.
 In contrast, the Cremona-convex hull (a subset of the frame-convex hull) is
 always 1-dimensional and forms part of the ruling of a quadratic cone with
 vertex $p$, except in these degenerate cases, when it coincides with the
 frame-convex hull. 
 We illustrate this when the lines have direction vectors $(1,2,3)$ and
 $(2,3,1)$. 
 \[
   \begin{picture}(250,140)
    \put(0,0){\includegraphics[width=5cm]{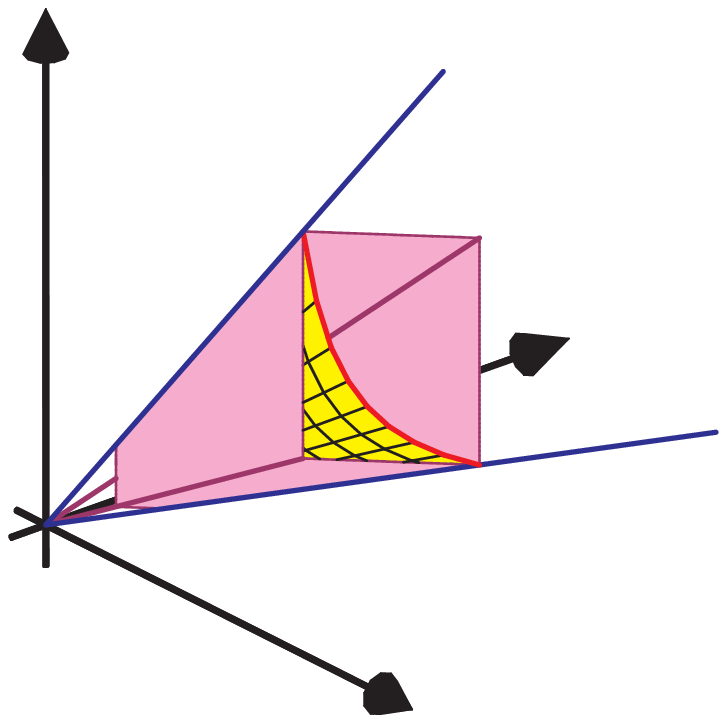}}
    \put(85,115){$\ell$}
    \put(130,40){$\ell'$}
    \put(75,17){$m$} \put(78,24){\vector(-1,1){26}}
    \put(48,126){$m'$} \put(54,123){\vector(1,-2){20}}
    \put(120,113){\vector(-2,-1){41}}
    \put(120,113){\vector(-1,-1){31}}
    \put(122,110){facets of frame-convex hull}
    \put(120,67){\vector(-4,-1){40}}
    \put(122,64){Cremona-convex hull}
  \end{picture}
\]

\section{The Cremona transformation for the Grassmannian of lines}\label{S:Cremona}

We describe some of the beautiful geometry behind the Cremona transformation,
which is induced by the Cremona coordinates~\eqref{E:Cremona_Coords}.
This allows us to describe the Cremona-affine hull of up to $d$ ascending lines
in $\R^d$.
A good reference for the vivid classical geometry that we use is Harris's
book~\cite{Ha92}.

\subsection{Geometry of the Cremona transformation}\label{S:linearizing}
The transformation between rectilinear and Cremona coordinates is
better understood in terms of homogeneous Pl\"ucker coordinates,
which are defined by
 \begin{equation}\label{E:Pluecker_Coordinates}
   p_{0i}\ =\ v_{i}\quad\textrm{and}\quad
   p_{ij}\ =\ x_{i}v_{j}-x_{j}v_{i}
  \qquad\textrm{for}\quad 1\leq i<j\leq d\,.
 \end{equation}
Let $\P$ be Pl\"ucker space, which has $\binom{d+1}{2}$ coordinates.
The {\it Cremona transformation} 
\[ 
  C\colon \P \to \P\quad {\rm is\, given \, by}\quad
   q\ =\ C(p)\,,
\]
 where
 \begin{equation}\label{E:Cremona_Transformation}
   q_{0i}\ =\ \frac{1}{p_{0i}}\quad\mbox{\rm and}\quad
   q_{ij}\ =\ \frac{p_{ij}}{p_{0i}p_{0j}}\,,
   \qquad\mbox{\rm for}\quad 1\leq i<j\leq d\, .
 \end{equation}
Let $\ell$ be a line whose direction vector has no vanishing
component.
If $(x,v)$ are rectilinear coordinates of $\ell$, $(y,w)$ are
Cremona coordinates of $\ell$ as in~\eqref{E:Cremona_Coords}, and
$q$ are its Cremona coordinates as in~\eqref{E:Cremona_Transformation}, then 
 \begin{equation}\label{E:new_Cremona}
   q_{0i}\ =\ w_i\qquad\textrm{and}\qquad q_{ij}\ =\ y_i-y_j\,,
   \qquad\textrm{for}\quad 1\leq i<j\leq d\, , 
 \end{equation}
which does not depend upon the choice of $y$.
This Cremona transformation is an isomorphism outside the
coordinate hyperplanes $p_{0i}=0$ for $i=1,\dotsc,d$.
It is also an involution, since applying it twice on points where it is
an isomorphism gives the identity.

\begin{rem}
 While this agrees with the notion of Cremona coordinates of~\eqref{E:Cremona_Coords}
 if no Pl\"ucker coordinate $p_{0i}$ vanishes, these two notions differ for lines
 where some of these coordinates vanish.
 (Recall that $(p_{01},\dotsc,p_{0d})$ is the direction vector of a
 line with Pl\"ucker coordinates $p_{ij}$). 
 That is because the Cremona tranformation is continuous where it is defined
 (which is described in Theorem~\ref{T:Cremona}),
 whereas the assignment of Cremona coordinates~\eqref{E:Cremona_Coords} is not 
 continuous (but defined everywhere). 
\end{rem}

The Grassmannian of lines $G(1,d)$ in Pl\"ucker space  is defined by
the Pl\"ucker relations
\[
   \begin{array}{rl}
    p_{0i}p_{jk}-p_{0j}p_{ik}+p_{0k}p_{ij}\ =\ 0\ &
      \textrm{for}\quad 1\leq i<j<k\leq d,\qquad \textrm{and}\\
    p_{ij}p_{kl}-p_{ik}p_{jl}+p_{il}p_{jk}\ =\ 0\ &\rule{0pt}{14pt}
      \textrm{for}\quad 1\leq i<j<k<l\leq d\,.
   \end{array}
\]
Let $1\leq i<j<k\leq d$.
Dividing the Pl\"ucker relation
\[
  {p_{0i}}p_{jk}-p_{0j}p_{ik}+p_{0k}p_{ij}\ =\ 0
\]
 by $p_{0i}p_{0j}p_{0k}$ reduces it to the linear relation
\[
  q_{jk}-q_{ik}+q_{ij}\ =\ 0\,.
\]
Let $1\leq i<j<k<l\leq d$.
Dividing the Pl\"ucker relation
\[
  p_{ij}p_{kl}-p_{ik}p_{jl}+p_{il}p_{jk}\ =\ 0
\]
by $p_{0i}p_{0j}p_{0k}p_{0l}$ reduces it to
 \begin{multline*}
 \quad q_{ij}q_{kl}-q_{ik}q_{jl}+q_{il}q_{jk}\ =\ \\
   q_{kl}(q_{ij}-q_{ik}+q_{jk})\ +\
    q_{ik}(q_{jk}-q_{jl}+q_{kl})\ +\
    q_{jk}(q_{ik}-q_{il}+q_{kl})\ =\ 0.\quad
 \end{multline*}
Therefore the Cremona transformation maps $G(1,d)$ into the linear
subspace $LG(d)$:
\[
   q_{jk}-q_{ik}+q_{ij}\ =\ 0\,,
   \qquad\textrm{for}\quad 1\leq i<j<k\leq d\, .
\]
We may use these linear relations to express any $q_{jk}$ with $1<j$ in terms of
the coordinates $q_{1i}$.
Since no equation defining $LG(d)$ involves any $q_{0i}$, these relations define
a linear subspace of dimension $2d-2$, which equals the dimension of the
Grassmannian.
Since the Cremona transformation is an isomorphism outside the
coordinate hyperplanes $p_{0i}=0$, we see that $G(1,d)$ is mapped birationally
onto the linear space $LG(d)$.

\begin{thm}\label{T:Cremona}
 The birational Cremona transformation $C$ on Pl\"ucker space is an involution
 that maps $G(1,d)$ to $LG(d)$ birationally.
 The indeterminancy locus of $C$ has codimension $3$ and consists of the linear
 subspaces
\[
  \begin{array}{rrl}
     L_{ijk}\ :&\quad p_{0i}\ =\ p_{0j}\ =\ p_{0k}\ =\ 0\ &
      \textrm{for}\quad 1\leq i<j<k\leq d,\qquad \textrm{and}\\
      L_{ij}\ :&\quad  p_{ij}\ =\ p_{0i}\ =\ p_{0j}\ =\ 0\rule{0pt}{14pt}
      &\textrm{for}\quad 1\leq i<j\leq d\,.
  \end{array}
\]
 The indeterminacy locus of the Cremona transformation restricted to
 $G(1,d)$ has $\binom{d}{2}+1$ components, all of codimension $2$. One
 component is the set of lines at infinity, while the other components
 consist of the sets of lines whose direction vectors have at least two
 coordinates equal to $0$.

 The indeterminancy locus of the Cremona transformation restricted to $LG(d)$
 has codimension $3$ and consists of the intersections of the subspaces
 $L_{ijk}$ and $L_{ij}$ with
 $LG(d)$.

 This indeterminacy locus in $LG(d)$ is the image under $C$
 of lines whose direction vectors have one coordinate equal to $0$.
\end{thm}

The set of lines in $\P^d$ that meet a linear subspace having codimension $a{+}1$ is
a {\em special Schubert cycle}.
This irreducible subvariety has codimension $a$ in $G(1,d)$.
More generally, given subspaces $L\subsetneq M$ of $\P^d$ where $L$ has
codimension $a{+}1$ and $M$ has codimension $b$, then the set of lines lying in
$M$ that meet $L$ is a Schubert cycle.
This has codimension $a{+}b$ in $G(1,d)$ and is special when $b=0$.\medskip 

\noindent{\it Proof.}
 We have already noted that Cremona transformation is an involution that restricts
 to a birational map from $G(1,d)$ to $LG(d)$.
 If we multiply the coordinates of the Cremona
 transformation~\eqref{E:Cremona_Transformation} by
 $p_{01}p_{02}\dotsb p_{0d}$, we see that it is also given by
 \begin{equation}\label{E:Cremona_polys}
     \begin{array}{rcll}
     q_{0i}&=& p_{01}\dotsb \widehat{p_{0i}}\dotsb p_{0d} &
      \textrm{for}\quad 1\leq i\leq d,\qquad \textrm{and}\\
     q_{ij}&=& p_{ij}p_{01}\dotsb \widehat{p_{0i}}\dotsb
                \widehat{p_{0j}}\dotsb p_{0d} \rule{0pt}{14pt}
      &\textrm{for}\quad 1\leq i<j\leq d\,,
  \end{array}
 \end{equation}
 and is therefore undefined if either
\begin{enumerate}
 \item[(i)] any three coordinates $p_{0i},p_{0j},p_{0k}$ vanish, or
 \item[(ii)] any two coordinates  $p_{0i},p_{0j}$ vanish, and we have
       $p_{ij}=0$.
\end{enumerate}
 But (i) defines the subspace $L_{ijk}$ and (ii) defines the subspace $L_{jk}$.

The restriction of these linear spaces to $G(1,d)$ is conveniently
described using special Schubert cycles.
The hyperplane section $p_{0i}=0$ is the set of lines that meet the
codimension 2 linear subspace $z_0=z_i=0$, which is a hyperplane at infinity.  
Similarly, the hyperplane section $p_{ij}=0$ is the set of lines that meet the
codimension $2$ linear subspace defined by $z_{i}=z_{j}=0$.  
Therefore $L_{ijk}\cap G(1,d)$ has two components:  the set of lines at
infinity, and the set of lines that 
meet the codimension $4$ subspace $z_{0}=z_{i}=z_{j}=z_{k}=0$.  
The first of these sets is a special Schubert
cycle of codimension $2$, while the other is a special Schubert cycle
of codimension $3$ in $G(1,d)$.  
Similarly, $L_{ij}\cap G(1,d)$ is the set of lines that
meet the codimension $3$ linear subspace $z_{0}=z_{i}=z_{j}=0$.
This is a special Schubert cycle of codimension $2$ that
contains the second component of the intersection  $L_{ijk}\cap G(1,d)$.

 Finally, the linear forms defining either $L_{ijk}$ or $L_{ij}$ are linearly
 independent of those defining $LG(d)$, therefore they define
 codimension $3$ linear
 subspaces of $LG(d)$ as well.  This indeterminacy locus is the image
 under $C$ of the coordinate hyperplane sections $p_{0i}=0$ of
 $G(1,d)$.
 The direction vector of such lines has some coordinate equal to $0$.
\QED

The restriction of the Cremona transformation to coordinate
hyperplanes $p_{0i}=0$, outside the indeterminacy locus, is a
contraction.  

\begin{prop}
 The coordinate hyperplane $p_{0i}=0$ meets the Grassmannian
 $G(1,d)$ in the Schubert cycle of lines that meet the codimension $2$ plane
 $z_0=z_{i}=0$ at infinity. It is contracted by $C$ to a line in $LG(d)$.
\end{prop}

\noindent{\em Proof.}
 The Pl\"ucker coordinates of a line that meets the codimension 2 plane 
 $z_0=z_i=0$ have $p_{0i}=0$, so its only nonzero Cremona coordinates are $q_{0i}$
 and $q_{1i}=\dotsb=q_{id}$. 
\QED

\subsection{Cremona-affine span of lines}\label{S:Cremona-affine}
Recall that if $f\colon X\, \to Y$ is a rational map with $Z\subset X$ an
irreducible subvariety that is not contained in the indeterminacy locus
$W\subset X$ of the map $f$, then the (proper) image of $Z$ is the closure
$\overline{\varphi(Z\setminus W)}$ in $Y$ of its set-theoretic image.
Given a set $\Gamma$ of lines, their Cremona coordinates span a linear subspace
$L(\Gamma)$ in $LG(d)$ whose image $G(\Gamma)$ in $G(1,d)$ under the Cremona
transformation is the {\em Cremona-affine span} of the original lines.
We use the Cremona map of Section~\ref{S:linearizing} to obtain a
precise description of such Cremona-affine spans.

Since the indeterminacy locus of the Cremona transformation has
codimension $3$ when restricted to $LG(d)$ we are able to describe
the subvarieties of the Grassmannian
that correspond to general linear spaces of lines of dimension up to $3$
in Cremona coordinates. A $2$- or $3$-dimensional subvariety $X$ of the
Grassmannian has a bidegree $(d_{1},d_{2})$ given by its intersection with general
Schubert cycles of codimension $2$ and $3$ respectively. In
codimension $2$, $d_{1}$ counts the number of lines in $X$ that meet a
general codimension $3$ linear space, while $d_{2}$ counts the
number of lines that lie in a general hyperplane.  In
codimension $3$, $d_{1}$ counts the number of lines in $X$ that meet a
general codimension $4$ linear space, while $d_{2}$ counts the
number of lines that lie in a hyperplane and meet a codimension $2$
linear subspace of this hyperplane.

A linear subspace in $LG(d)$ is in {\it general position} if it meets the
indeterminacy locus of $C$ properly, that is, in a set of codimension $3$.  
The Cremona coordinates $q_{0i}$ of a line are the
inverses of the components of its direction vector.
We call the vector $(q_{01},\dotsc,q_{0d})$ the 
{\it Cremona direction} of the line.   A set
of $n \leq d$ lines are {\it independent} if their Cremona
directions are affinely independent.
If their affine span is in addition in general position in $LG(d)$, then they
are in general position.

Note that a set of independent lines have affinely independent images
in $LG(d)$, while the converse is not true.  The converse fails for
instance when two lines have the same direction.

\begin{thm} Let $\Gamma$ be a set of $n$ independent lines in general position
  in $\R^d$.
  Then their Cremona-affine span $G(\Gamma)$ is a rational $(d{-}1)$-fold in
  $G(1,d)$.  

  \begin{enumerate}
       \item  If $n=2$, then $G(\Gamma)$ is a rational normal curve of degree
         $d{-}1$. 

       \item If $n=3$, then $G(\Gamma)$ is a Veronese surface of bidegree
          $(\binom{d-1}{2},\binom{d}{2})$ and degree $(d-1)^2$.

       \item If $n=4$, then $G(\Gamma)$ is a rational threefold of bidegree
       $(\binom{d-1}{3}, 2\binom{d}{3})$ and degree
       $(d{-}1)^3-\binom{d}{3}-\binom{d}{2}$.
   \end{enumerate}
   \end{thm}

\noindent{\it Proof.}
 Clearly these varieties are rational.
  For the numerical results, we compute the bidegree and degree of $G(\Gamma)$.
 The Cremona transformation as expressed in~\eqref{E:Cremona_polys}  is
 defined by polynomials of degree $d{-}1$.
    The indeterminacy locus has codimension $3$, so for $k=2$ and
    $k=3$ the degree of $G(\Gamma)$ is $(d{-}1)^{k-1}$.  For $k=3$ the
    indeterminacy locus consists of the points of intersection with
    $L_{ijk}$ and $L_{ij}$, which add up to $\binom{d}{3}+
    \binom{d}{2}$ points, from which the degree of $G(\Gamma)$ follows.  The
    bidegree is determined by the degree of $G(\Gamma)$  and the degree
    of the corresponding union of lines in $\R^d$.  The latter can be
    computed by considering a codimension $k{-}1$ subspace of the
    hyperplane at infinity, and counting how many lines in the family
    meet this subspace.  This  is determined by the direction
    vector of the line, which consists of its Pl\"ucker coordinates
 $p_{0i}$~\eqref{E:Pluecker_Coordinates}.
    On
    these coordinates the Cremona transformation is the ordinary
    Cremona transformation defined by inverting all the
    coordinates.  The degree of the restriction of an ordinary Cremona
    transformation on $\P^m$ to a general $n$-dimensional linear
    subspace is $\binom{m}{n}$. Therefore the degree of the image of a
    $(k{-}1)$-dimensional linear subspace at infinity is $\binom{d-1}{k-1}$.
\QED

We describe the Cremona-affine hull
of a finite set of lines whose direction vectors have nonzero
coordinates, as a subset of $\R^d$.
Consider the homogeneous Pl\"ucker coordinates $p_{ij}$
for an ascending line $\ell$.  We may choose two points on the line
with coordinates
\[
\begin{pmatrix}
    1 & 0 & -p_{12} & \ldots & -p_{1d}  \\
    0 & 1 & p_{02} & \ldots & p_{0d}
\end{pmatrix}
\]
We have scaled the Pl\"ucker coordinates so that $p_{01}=1$.
This is possible by the assumption that the coordinates $p_{0i}$ are
nonzero.
The linear forms vanishing on $\ell$ have a basis
 \begin{eqnarray*}
   b_{i}(\ell)&=&\frac{p_{1i}}{p_{0i}}z_{0}-z_{1}+ \frac{1}{p_{0i}}z_{i}\\
    &=& q_{1i}z_0-z_1+q_{0i}z_i\,,
   \qquad\textrm{for}\quad i=2,\dotsc,d\, .
 \end{eqnarray*}
These are Santal\'o's equations~\eqref{E:Cremona_revelaled},
expressed in the Cremona coordinates~\eqref{E:new_Cremona}.

Given a set of $n$ lines $\Gamma=\{\ell_{1},\ldots,\ell_{n}\}$,
let $M(\Gamma)$ be the $n\times (d{-}1)$-matrix 
of linear forms with entries $b_{i}(\ell_{j})$.
As in Remark~\ref{R:Santalo_proof}, a linear combination of the rows of
$M(\Gamma)$
gives a vector of linear forms whose coefficients are
the same linear combination of the Cremona coordinates of the $\ell_i$.
This identifies the row space of $M(\Gamma)$ with the affine span of
the Cremona coordinates of the lines $\ell_{i}$.
The matrix $M(\Gamma)$ provides a determinantal description of the
variety of lines $S(\Gamma)\subset \R^d$ that belong to the linear span of
$\Gamma$ in Cremona
coordinates.

\begin{thm}\label{T:Nonzero}
   Let $\Gamma$ be a set of $n$  independent lines in $\R^d$ whose direction
   vectors have only nonzero coordinates.  
   Let $L(\Gamma)\subset LG(d)$ be the linear span of their Cremona coordinates,
   $G(\Gamma)$ their Cremona-affine span, and $S(\Gamma)\subset\R^d$ the
   union of the lines in $G(\Gamma)$.
 \begin{enumerate}
  \item The set of lines in the family $G(\Gamma)$ that pass through a
     point in $\R^d$ is linear in $L(\Gamma)$.
  \item When $n\leq d{-}1$, $S(\Gamma)$ is a $n$-fold scroll of degree
    $\binom{d-1}{n-1}$ defined by the maximal minors of $M(\Gamma)$.
  \item The matrix $M(\Gamma)$ has rank at most $r$ precisely at the
    points in $\R^d$ that are contained in a $(n{-}r{-}1)$-dimensional set of 
    lines in the family.
  \item
    If $n\leq d$, there is at most one line in the family
    in each direction with only nonzero coordinates.
  \item  For $n = d$ the set of points in $\R^d$ that lie on more than one
    line in the family form a codimension $2$ subvariety $Z(\Gamma)$ in $\R^d$ of
    degree $\binom{d}{2}$.  In particular, $G(\Gamma)$ is
    the set of $(d{-}1)$-secant lines to $Z(\Gamma)$.
  \end{enumerate}
 \end{thm}

Part (2) for $n=2$ and $d=3$ and Part (5) for $d=3$ are illustrated in
  Figure~\ref{F:Cremona_Interval}.
  Notice that the degree of the union of lines in $\R^d$ of the
  family $G(\Gamma)$ coincides
  with the number of lines in the family that meet a general
  codimension $n$ linear subspace.\smallskip

\noindent{\em Proof.}
 The set of lines in the family that pass though a point
are cut out by linear combinations of rows in $M(\Gamma)$ that vanish when
evaluated at that point, so this set is a linear subspace of $L(\Gamma)$.  
If the rank of the matrix $M(\Gamma)$ at $p$ is $r<n$, then the dimension of
this set of lines is $n{-}r{-}1$. 
The union of all lines in the family is precisely the $n$-dimensional subvariety
of $\R^d$ where the matrix has rank at most $n{-}1$.  
Its  degree is given by the Thom-Porteus formula~\cite[p.~254]{Fu}.
When $n=d$, the matrix $M(\Gamma)$ has
dimensions $n\times (n-1)$.  It therefore has rank at most $n-2$ on
a $(d{-}2)$-dimensional  subvariety  $Z(\Gamma)$ of degree $\binom{d}{2}$ in $\R^d$,
 defined by the maximal minors
of $M(\Gamma)$.  For each line in the family, this intersection with
$Z(\Gamma)$ is defined by a unique $(d{-}1)$-dimensional minor, so the
line is $(d{-}1)$-secant to $Z(\Gamma)$.  On the other hand, let $L$ be a
$(d{-}1)$-secant line to $Z(\Gamma)$.  Then the subspace of maximal
minors of $M(\Gamma)$ that vanish on $\ell$  has codimension $1$.  But
the codimension $1$ subspaces of rows of the matrix correspond
naturally to this space
of minors, so the codimension $1$ subspaces
corresponding to minors that vanish on $\ell$ have a unique row in
common.  This row must define the line $\ell$.
\QED

Let $n< d$.
Note that the matrix $M(\Gamma)$ has rank 1 at infinity along 
the $i$-th coordinate axis (where $z_{j}=0$ for $j\neq i$).
Therefore there is a $(n{-}1)\times(d{-}1)$-submatrix
$M_{i}(\Gamma)$ of $M(\Gamma)$ that has rank $0$
at infinity along the $i$-th coordinate axis.

\begin{cor}  
   Let $\Gamma$ be a set of $n$ independent lines whose direction vectors have
   only nonzero coordinates in $\R^d$, where $1<n<d$. 
   The lines in the family $G(\Gamma)$ that are parallel to the $i$-th
   coordinate axis form a 
   $(n{-}2)$-dimensional family $G_{i}(\Gamma)\subset G(1,d)$ which is the
   preimage under $C$ of the codimension $1$ linear subspace
   $L_{i}(\Gamma)=\{q_{0i}=0\}\cap L(\Gamma)$.
   The corresponding subvariety $S_{i}(\Gamma)\subset S(\Gamma)$  is an
   $(n{-}1)$-dimensional cylinder over an $(n{-}2)$-dimensional
   subvariety of the coordinate hyperplane $\{z_{i}=0\}$ defined by the maximal
   minors of the $(n{-}1)\times(d{-}1)$-matrix $M_{i}(\Gamma)$.  
   In particular,  $S_{i}(\Gamma)$ is a cylinder over a
   rational determinantal variety of degree  $\binom{d-1}{n-2}$.
\end{cor}

We characterize some sets $\Gamma$ of independent lines
and some general linear subspaces $L({\Gamma})$.

\begin{prop}\label{P:indep}
    Two lines whose direction vectors have only nonzero coordinates are
independent if their direction vectors are distinct. Three lines whose
direction vectors have only nonzero coordinates are independent if and
only if there is no surface of minimal degree $d{-}1$ that contains
all three and passes through the coordinate directions at infinity.
\end{prop}

\noindent{\em Proof.} 
 The first statement is the definition of independent.
 For the second statement, notice
 that the construction of the matrix $M(\Gamma)$ defines
 the unique scroll of minimal degree that contains two lines whose direction
 vectors have only nonzero
 coordinates and that passes through the coordinate directions at infinity.
\QED

Given two lines in $\R^d$ whose direction vectors have only nonzero
coordinates, then the scroll of degree $d{-}1$ that is defined by their
Cremona-affine span contains one line parallel to each coordinate axis.
This line has a direct synthetic description:
Consider two lines $\ell_{1}, \ell_{2}$, and the coordinate simplex
$\Delta$ at infinity.  
For each codimension $2$ subsimplex $\Delta'\subset\Delta$, consider
the codimension $2$ subspace $P_{\Delta'}$ through $\Delta'$ that meets both lines
$\ell_{1}, \ell_{2}$.  
There is a unique line through each vertex of $\Delta$ 
that meets all the subspaces $P_{\Delta'}$.  This is the
line of the scroll parallel to the chosen coordinate axis.

Together with the original  lines $\ell_{1}, \ell_{2}$ we count $d+2$ lines that
meet $\binom{d}{2}$ codimension $2$ linear
spaces  $P_{\Delta'}$.  The set of lines that meet a codimension $2$
linear space form a hyperplane section of $G(1,d)$.  In general, the hyperplanes
defined
by the different $P_{\Delta'}$
are linearly independent, and therefore define a $\P^{d-1}$ inside
Pl\"ucker space.    The $d+2$ lines are represented by  $d+2$ points in this
space.
By Castelnuovo's lemma~\cite[p.~530]{GH78}, there is a unique rational normal
curve through these $d+2$ points.  
This curve defines the scroll $S(\ell_1,\ell_2)$ of Theorem~\ref{T:Nonzero}.

It is straightforward to verify the following properties of
the scroll  $S(\ell_1,\ell_2)$.

\begin{prop}\label{P:degenerate}
  Let $\ell_1 $ and $ \ell_2$ be two lines whose direction vectors 
  have only nonzero coordinates. If they are disjoint,
  then their span $S(\ell_1,\ell_2)$ is a smooth scroll.
  If they meet at a finite point, then their span is a cone over a
  rational normal curve.  
\end{prop}

Example~\ref{Ex:lines_meet} illustrated the case of
Proposition~\ref{P:degenerate} in $\mathbb{R}^3$ where the two lines 
have non-vanishing coordinate directions and they meet. 

\begin{prop}\label{prop:general}
 Let $\ell_{1}, \ell_{2}$, and $\ell$ be three lines in
 $\R^d$ whose
 direction vectors have only nonzero coordinates.   Assume that the two 
 lines $\ell_{1}$ and $\ell_{2}$ are
independent and  have
direction vectors whose affine span does not meet any codimension $2$
coordinate plane.  Then the following are
 equivalent:
\begin{enumerate}
    \item $\ell$ meets every codimension $2$ plane parallel to two
    coordinate hyperplanes that also meet $\ell_{1}$ and $\ell_{2}$.
    \item The projections of the three lines into any coordinate plane
    have a common point.
    \item $\ell$ is a ruling in the unique scroll of minimal degree
    that contains $\ell_{1}$ and $\ell_{2}$ and passes through the
    coordinate directions at infinity.
    \end{enumerate}
\end{prop}

\noindent{\em Proof.} 
Since the span of the direction vectors of $\ell_1$ and $\ell_2$ does
not intersect any codimension $2$ coordinate plane at infinity, the
projection of the two lines to any coordinate plane meet in a unique
point. 
The preimage of this point is a codimension $2$ plane parallel to
two coordinate hyperplanes that meet both lines, and any such codimension
$2$ plane arises this way. Therefore (1) and (2) are equivalent.  

Let $H_{ij}$ be the unique codimension 
$2$ plane that intersects both $\ell_{1}$ and $\ell_{2}$ and is parallel to all 
coordinate axes except the $i$- and the $j$-axis.    

For (3), note that by Propositions~\ref{P:indep} and~\ref{P:degenerate} the
 scroll $S{(\ell_1,\ell_2)}$ is smooth 
 and uniquely determined by the two lines and the coordinate directions at infinity.
 The intersection $C_{ij}=H_{ij}\cap S(\ell_1,\ell_2)$ contains a point at
 infinity on $d{-}2$  
 coordinate axes in addition to one point on each of  the lines  $\ell_1$ and $\ell_2$.  
 Altogether we have counted $d$ points, while the scroll has degree only $d{-}1$.  
 Therefore the intersection $C_{ij}$ contains a curve which furthermore must
 intersect every line   
 in the ruling of $S{(\ell_1,\ell_2)}$ defined by $\ell_1$ and $\ell_2$. 
 In particular, $\ell$ must intersect $H_{ij}$ if it belongs to the scroll
 $S{(\ell_1,\ell_2)}$. 
 
  On the other hand, the Schubert cycle of lines meeting $H_{ij}$
forms a hyperplane section $\Lambda_{ij}\cap G(1,d)$ of $G(1,d)$ in
$\mathbb{P}$.   
The image of $\Lambda_{ij}$ under the Cremona transformation is a hyperplane
$\Lambda'_{ij}$ that meets $LG(d)$ in a hyperplane $P_{ij}$ of $L(G)$.   
If $\ell$  meets each $H_{ij}$, its Cremona coordinates must lie in
every $P_{ij}$.  Parallel to each coordinate 
axis there is exactly one line that meets each $H_{ij}$, so 
the intersection of the hyperplanes $P_{ij}$ is exactly the line spanned by the
Cremona coordinates of $\ell_{1}$ and $\ell_{2}$.  
Therefore the Cremona coordinates of $\ell$ belong to this line, and $\ell$
belongs to $S{(\ell_1,\ell_2)}$. 
\QED

We will use this to show that the frame-convex hull of two general ascending
lines equals their Cremona-convex hull.

First note that these two notions of convexity are preserved by orthogonal
projections along coordinate directions.
Indeed, suppose that $\ell$ is a (necessarily) ascending line lying in the
frame-convex hull of two ascending lines $\ell_1$ and $\ell_2$ and let $\pi$ be
an orthogonal projection to a coordinate $d'$-plane.
Then $\pi(\ell)$ lies in the frame-convex hull of $\pi(\ell_1)$ and
$\pi(\ell_2)$.
The same is true for Cremona convexity.

Suppose that $\ell_1$ and $\ell_2$ are ascending lines with different direction
vectors.
Then their Cremona-affine hull $G(\ell_1,\ell_2)$ consists of a rational curve
of lines with different directions, and the lines $\ell_1$ and $\ell_2$ divide
this curve (topologically a circle) into two intervals, one of which consists
of ascending lines.

\begin{cor}\label{cor:general}
 Suppose that $\ell_1$ and $\ell_2$ are ascending lines with different
 direction vectors.
 Then their Cremona-convex hull is the interval they define on
 $G(\ell_1,\ell_2)$ consisting of ascending lines. 
 If $d=2$ or if the two lines are disjoint and have
direction vectors whose affine span does not meet any codimension $2$
coordinate plane, then this is also their
 frame-convex hull. 
\end{cor}

\noindent{\it Proof.}
 The first statement is straightforward:
 Since $\ell_1$ and $\ell_2$ are ascending, every line in their Cremona-convex
 hull is also ascending.
 But this Cremona-convex hull is one of the intervals they define on
 $G(\ell_1,\ell_2)$. 

 Suppose now that either $\ell_1$ and $\ell_2$ are skew and have
direction vectors whose affine span does not meet any codimension $2$
coordinate plane or else $d=2$.
 By assumption, their images under any orthogonal projection $\pi$ to a coordinate $2$-plane
 meet in a unique point $p$, and the inverse image $\pi^{-1}(p)$ is a codimension $2$ plane $H$ that is parallel to
 the codimension $2$ coordinate plane $\pi^{-1}(0)$.
 Let $B(\subset H)$ be the smallest box containing the points $\ell_1\cap H$ and
 $\ell_2\cap H$ (necessarily as vertices).
 Thus the frame-convex hull of $\ell_1$ and $\ell_2$ is a subset of the set of
 lines meeting each such coordinate-parallel codimension 2 plane $H$, but this
 larger set is the Cremona-affine hull of $\ell_1$ and $\ell_2$, by
 Proposition~\ref{prop:general}.

 Choose a different projection $\pi'$ to a coordinate 2-plane in which the box
 $B$ projects to an interval of positive length.
 The endpoints of this interval are necessarily the projections of the points 
 $\ell_1\cap H$ and $\ell_2\cap H$, and we have 
 $\pi'(\ell_1)\neq\pi'(\ell_2)$.
 Since $\pi'$ is 1-1 on lines with no vanishing coordinate directions, in
 particular on ascending lines, it maps lines in the frame-convex hull of $\ell_1$
 and $\ell_2$ to lines in the frame-convex hull of $\pi'(\ell_1)$ and
 $\pi'(\ell_2)$. 
 This is the interval of lines in that coordinate plane through
 $p:=\pi'(\ell_1)\cap\pi'(\ell_2)$ consisting of ascending lines, and so is also
 the image of the Cremona-convex hull of $\ell_1$ and $\ell_2$.
 Since $\pi'$ is 1-1 on ascending lines, this shows that the frame-convex hull of
 $\ell_1$ and $\ell_2$ equals their Cremona-convex hull.
\QED

\section{Convexity structures for affine subspaces}\label{S:Frame}

 We discuss convexity structures for general affine subspaces of
 $\R^d$. 
 We first describe Santal\'o's Helly-type theorem for hyperplanes, and then give
 a Helly-type result for codimension 2 linear subspaces via their Cremona
 coordinates. 
 Next, we discuss frame convexity for lines with different sign patterns, and
 finally show that  frame convexity for $k$-flats in $\R^d$ with $1<k<d{-}1$ has
 a completely different character than for points, lines, or hyperplanes.\smallskip
  
 Santal\'o also proved the following Helly-type theorem for hyperplane
 transversals.

\begin{thm}\label{T:18}
 If every $2^{d-1}(d+1)$ (or fewer) members of a collection of boxes in
 $\R^d$ admit a hyperplane transversal, then there is a hyperplane transversal
 to all the boxes.
\end{thm}

 As with Theorem~\ref{T:Santalo_Lines}, this will be implied by a local version.
 A hyperplane $H$ is {\it ascending} if it has equation
\[
    z\cdot v\ =\ \sum_{i=1}^d z_i v_i\ =\ c\,,
\]
 where the perpendicular vector $v=(v_1,\dotsc,v_d)$ is ascending.
 An ascending hyperplane meets a box with minimum point $a$ and maximum point
 $b$ if and only if
\[
   a\cdot v\ \leq\  c\ \leq\ b\cdot v\,.
\]
 Thus the set of ascending hyperplanes meeting this box is convex, in the
 coordinates given by the coefficients of their defining equations.
 Theorem~\ref{T:18} follows immediately from Helly's theorem, 
 in the same way as Theorem~\ref{T:Santalo_Lines}.

\begin{rem}
 The convexity for both ascending lines and ascending hyperplanes in $\R^d$ came
 from natural coordinates that identified them as convex subsets of
 $\R^{2d-2}$ and $\R^d$, respectively.
 For hyperplanes, this is quite natural, as the set of hyperplanes in $\P^d$ is
 just the dual $\P^d$.
 For lines, the Cremona transformation linearizes the Grassmannian of lines, and
 our notion of Cremona convexity is pulled back from the affine structure on
 this linearized Grassmannian.
 For $1<k<d-2$, there is no linearizing Cremona transformation for the
 Grassmannian of $k$-flats in $\P^d$, which partially explains our inability to
 extend the results for lines and hyperplanes to $k$-flats for arbitrary $k$.

 The other ingredient in these results is that the natural convex structure from
 the linearizations has the nice geometric interpretation of frame convexity.
 Since the set of $(d{-}2)$-flats is isomorphic to the space of lines, there is
 a linearizing Cremona transformation and a notion of convexity for
 $(d{-}2)$-flats.
 It does not, however have as nice a geometric interpretation as
 frame convexity.
 We describe it by dualizing the notion of frame convexity.

 The dual of a $k$-flat in $\R^d$ is a $(d{-}1{-}k)$-flat in $\RP^d$ that does
 not meet the origin.
 We restrict ourselves to flats in $\RP^d$ that do not meet the origin.
 Every hyperplane
 (at finite distance)
 not containing the origin has a unique equation of the form
\[
   a\cdot z\ =\ \sum_i a_i z_i\ =\ 1\,,
\]
 where $z_1,\dotsc,z_d$ are the coordinates of $\R^d$.
 The coefficients $a=(a_1,\dotsc,a_d)$ give
 coordinates for hyperplanes not containing the origin, with $(0,\dotsc,0)$
 giving the hyperplane at infinity.
 The set of such hyperplanes is partially ordered by componentwise comparison of
 these coordinates.
 This has a geometric interpretation.

 Order the non-zero points on a coordinate axis in $\P^d$ so that the positive
 numbers with their usual order precede the point at infinity, which precedes
 the negative numbers with their usual order.
 Given two hyperplanes $A$ and $B$ not containing the origin, we have that
 $A\leq B$ if and only if for each coordinate axis $\ell$,
 $\ell\cap A$ precedes $\ell\cap B$ in this order.
 The {\it $*$-box} between two hyperplanes $A\leq B$ not containing the origin
 is the set of those hyperplanes $H$ (not containing the origin) that satisfy
 $A\leq H\leq B$.

 A $(d{-}2)$-flat $K$ not meeting the origin is {\it ascending} if its span
 with the origin is an ascending hyperplane.
 That is, $K$ is the set of points $z$ in $\R^d$ defined by the equations
\[
   z\cdot x\ =\ 1,\qquad\mbox{and}\qquad z\cdot v\ =\ 0\,,
\]
 where each coordinate of $v$ is non-negative and $v\neq 0$.
 This pair $(x,v)$ gives {\it rectilinear coordinates} for $K$, and we obtain
its
 Cremona coordinates via the transform~\eqref{E:Cremona_Coords}.

 The hyperplanes containing $K$ that do not contain the origin have the
 form
\[
   H_t\ \colon\ \{z\in\R^d\mid z\cdot( x+tv)\ =\ 1\}\,,
   \qquad\mbox{for}\ t\in\R\,.
\]
 An ascending $(d{-}2)$-flat {\it $*$-transverse} to a $*$-box defined by
hyperplanes
 $A\leq B$ is an ascending 2-plane lying on a hyperplane $H$ in the $*$-box.

 Since these definitions are just a translation of those for lines to $(d{-}2)$-flats
 via duality, we have the following results.

\begin{thm}
  The set of Cremona coordinates of ascending $(d{-}2)$-flats that are
  $*$-transverse to a given $*$-box is convex.
\end{thm}

\begin{cor}
  If we have a collection of $*$-boxes in $\R^d$ such that every $2d-1$ of
  them have a $*$-transversal, then they all do.
\end{cor}

 Note that $*$-transversality makes sense for any $(d{-}2)$-flat that does not
 contain the origin.
 Here is Santal\'o's theorem in this context.

\begin{thm}
 If every $2^{d-1}(2d-1)$ (or fewer) members of a collection of $*$-boxes in
 $\R^d$ admit a $*$-transversal by a $(d{-}2)$-flat, then there is a
 $(d{-}2)$-flat $*$-transversal to all the $*$-boxes.
\end{thm}

\end{rem}
\begin{ex}
 While the frame-convex hull of two general ascending lines has dimension $1$,
 if the lines have different sign pattern, then their frame-convex hull is just
 the two lines again, as the figure below shows for the lines $\ell_1$ and
 $\ell_2$, and the three axis-parallel boxes (actually segments) $B_1$,
 $B_2$, and $B_3$.
\[
 \begin{picture}(110,110)
   \put(0,0){\includegraphics[width=4cm]{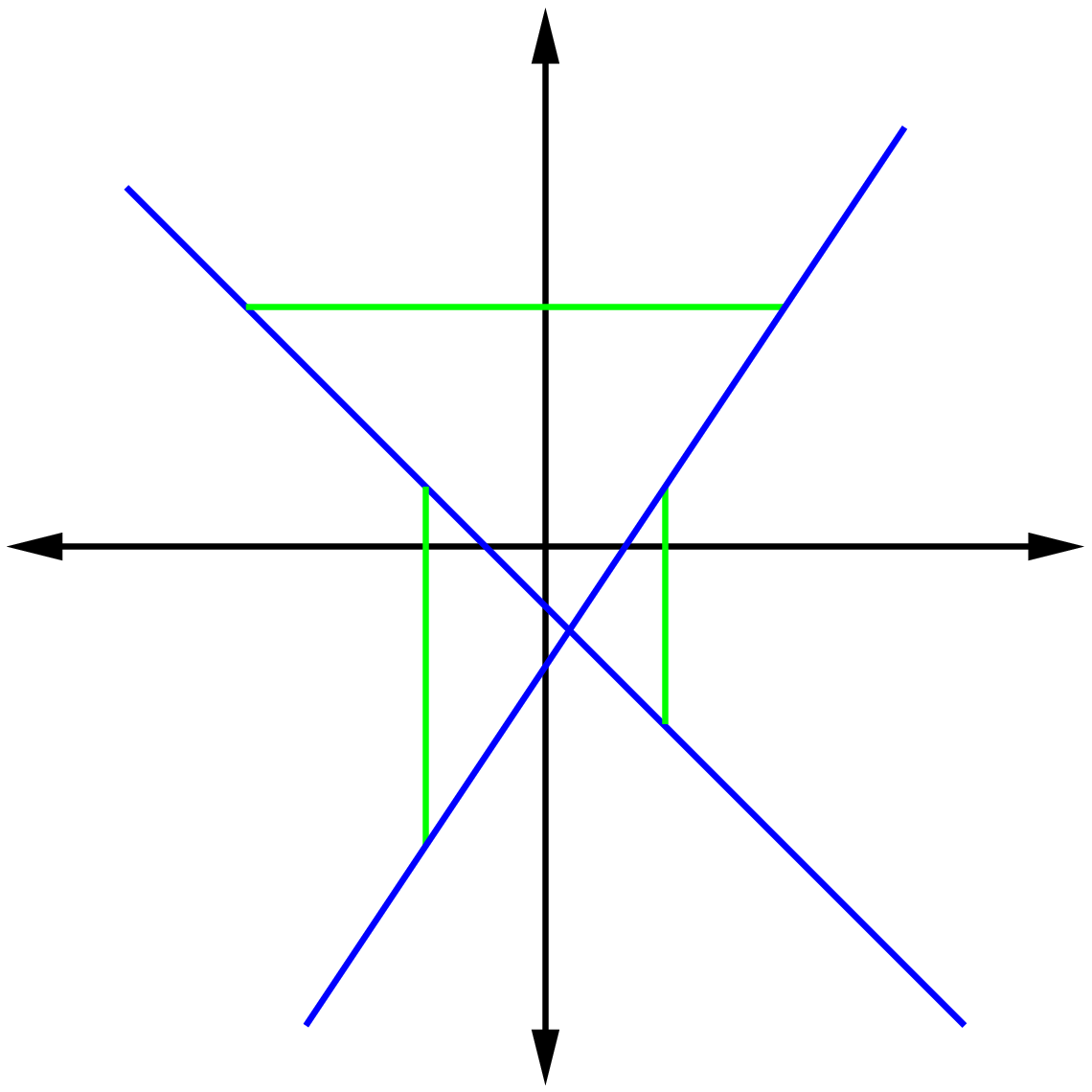}}
    \put(80,65){$\ell_1$}    \put(85,25){$\ell_2$}
    \put(30,40){$B_1$}       \put(35,85){$B_2$}         \put(70,45){$B_3$}
    \put(60,10){$y$}         \put(10,50){$x$}
  \end{picture}
\]
\end{ex}

 When $k\neq 1, d-1$, frame convexity for $k$-flats in $\R^d$
 has a completely different character than when $k$ is 1 or $d{-}1$.
 For example, the frame-convex hull of two 2-flats in either $\R^4$ or $\R^5$ is
 just the two original 2-flats.

\begin{prop}
 The frame-convex hull of any two $2$-flats in general position in either $\R^4$ 
 or $\R^5$ consists solely of the original $2$-flats.
\end{prop}

\noindent{\it Proof.}
 We give an algebraic relaxation of the frame-convex hull
 whose only solutions are the original 2-flats.
 A $k$-flat in $\R^d$ corresponds to a $(k{+}1)$-dimensional linear subspace of
 $\R^{d+1}$, which we represent as the row space of a $(k+1)\times (d+1)$
matrix,
 and consider as a $k$-flat in projective space $\RP^d$.
 Consider the two 2-flats in $\RP^5$
\[
     \left[\begin{array}{rrrrrr}
           1& 1& 1& 1& 1& 1  \\
           1&-1&-1& -1& 1& 0 \\
           0&-1& 1& 2& 2& 0  \end{array}\right]
      \qquad\mbox{and}\qquad
     \left[\begin{array}{rrrrrr}
          -1&-1& -1&-1& -1& 1 \\
          -1& 1& -1& 1& -1& 0 \\
           0&-1& 2& 1& 2& 0    \end{array}\right]\ .
\]
 For each coordinate $x_i$, there is a line parallel to the $x_i$-axis meeting
 both 2-planes.
 For each, we give the point in $\R\P^5$ on the line with $x_i=0$.
\[
   \begin{array}{rcl}
     \mbox{$x_1$-axis}&:&(0,-3,13,21,13,1)\\
     \mbox{$x_2$-axis}&:&(3,0,-5,-9,-5,1)\\
     \mbox{$x_3$-axis}&:&(-1,1,0,1,-1,0)\\
     \mbox{$x_4$-axis}&:&(1,1,-3,0,-3,0)\\
     \mbox{$x_5$-axis}&:&(9,1,-15,-23,0,1)
   \end{array}
\]
 Among the boxes meeting our two 2-flats are axis-parallel
 segments lying along these axis-parallel transversals.
 Thus the set $X$ of 2-flats meeting these axis-parallel transversals contains
the
 frame-convex hull of our original 2-flats.
 We consider $X$ to be an algebraic relaxation of the frame-convex hull.

 Formulating and solving the equations for a general 2-plane in $\P^5$
 to meet each of these five lines shows that $X$ consists of our original two
 2-planes.
 This is not unexpected: the space of 2-planes in $\P^5$ has dimension $9$, and
 the set of 2-planes meeting a line has codimension $2$, so there would be
 be no solutions to such a geometric problem, if the lines were in general
 position.
 Since this algebraic relaxation for a particular choice of two 2-flats has only
 the two original 2-flats as solutions, this will be the case for any two
 general 2-flats in $\R^5$.

 A different procedure shows that the solutions to an algebraic relaxation to
 the frame-convex hull of two 2-flats in $\R^4$ consist of only those 2-flats.
 If we choose a third 2-flat, then there are four axis-parallel lines meeting
 the three.
 Choosing yet another 2-flat gives four more axis-parallel lines, and an
 algebraic relaxation of the frame-convex hull of the original two 2-flats is
 the set of 2-flats meeting all 8 lines.
 Solving this in a specific instance gives the two original 2-flats,
 which likewise implies that the frame-convex hull of two general 2-flats in
 $\R^4$ consists only of the original 2-flats.
\QED

\def\cprime{$'$}
\providecommand{\bysame}{\leavevmode\hbox to3em{\hrulefill}\thinspace}
\providecommand{\MR}{\relax\ifhmode\unskip\space\fi MR }
\providecommand{\MRhref}[2]{%
  \href{http://www.ams.org/mathscinet-getitem?mr=#1}{#2}
}
\providecommand{\href}[2]{#2}

\end{document}